\documentclass{article}%
\usepackage{epsfig}
\usepackage{enumerate,graphicx}
\usepackage{amsfonts}
\usepackage{graphicx}
\usepackage{amsmath,amsthm,amssymb}
\usepackage{amssymb}%
\DeclareGraphicsExtensions{.pdf, .jpg, .png , .bmp, .ps}
\usepackage{tikz}
\providecommand{\U}[1]{\protect\rule{.1in}{.1in}}
\newtheorem{thm}{Theorem}[section]

\newtheorem{prop}{Proposition}[section]
\newtheorem{cor}{Corollary}[section]
\newtheorem{Defi}{Definition}[section]

\theoremstyle{remark}
\newtheorem{rem}[thm]{\bf Remark}

\newcommand\I{{\bf 1}}

\newcommand\bkE{{\mathbb {E}}}

\DeclareMathOperator{\Var}{V}

\def\1{{{\mbox{${\rm{1\negthinspace\negthinspace I}}$}}}}

\newcommand{\eref}[1]{(\ref{#1})}

\newcommand\beq{\begin{equation}}
\newcommand\eeq{\end{equation}}

\begin{document}

\title{Weak convergence of  the empirical process of intermittent maps in ${\mathbb L}^2$ under  long-range dependence.}

\author{J\'er\^ome Dedecker\footnote{Universit\'e Paris Descartes,
Sorbonne Paris Cit\'e,
Laboratoire MAP5
and CNRS UMR 8145, 75016 Paris, France.},
Herold Dehling\footnote{Faculty of Mathematics, Ruhr-University Bochum, 44780 Bochum, Germany.}  and
Murad S. Taqqu\footnote{Department of Mathematics,
Boston University,
Boston MA 02215, USA}}

\maketitle

\abstract{We study the behavior of the empirical distribution function of  iterates of intermittent maps in the Hilbert space
of square inegrable functions with respect to Lebesgue
 measure. In the long-range dependent case, we prove that
 the empirical distribution function, suitably normalized, converges
 to a degenerate stable process, and we give the corresponding almost sure result. We apply the results to
the convergence of the Wasserstein distance between the empirical measure and the
 invariant measure. We also apply it to obtain the asymptotic distribution of the corresponding Cram\'er-von-Mises statistic.}

\medskip


\medskip


\section{Introduction and main results.}\label{intro}
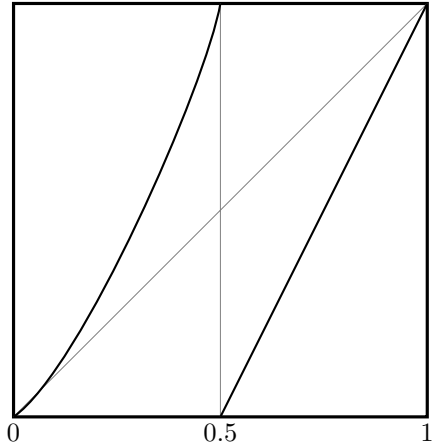
\begin{figure}[htb]
\label{fig:pm_graph}
\centering
  \begin{tikzpicture}[scale= 1.1]
  \draw[very thick] (0,0) rectangle (5,5);
  \draw[gray, very thin]
      (2.5,0) -- +(0, 5)
      (0,0)   -- (5,5);
  \draw[thick]
      (0,0) .. controls +(35:1) and +(-100:1) .. (2.5, 5)
      (2.5, 0) -- (5, 5);
\foreach \x/\ytext in {0/$0$, 1.2/$$, 2.5/$0.5$, 4.2/$$, 5/$1$}
\node[above] at (\x, -0.4) {\ytext};
  \end{tikzpicture}
\caption{Graph of intermittent map $T_\gamma: [0,1]\rightarrow [0,1]$}
\end{figure}

For  $\gamma$ in $]0, 1[$, we consider the intermittent map $T_\gamma$
(or simply $T$) from $[0, 1]$ to
$[0, 1]$, introduced   by Liverani, Saussol and Vaienti (1999):
$$
   T_\gamma(x)=
  \begin{cases}
  x(1+ 2^\gamma x^\gamma) \quad  \text{ if $x \in [0, 1/2[$}\\
  2x-1 \quad \quad \quad \ \  \text{if $x \in [1/2, 1]$;}
  \end{cases}
$$
see Figure~\ref{fig:pm_graph} for the graph of $T_\gamma$.

This kind of maps are known to exhibit a transition  from a stable periodic behavior to
a chaotic one, as described in  Pomeau and Manneville (1980).
Concerning the existence of $T_\gamma$-invariant probability measures which are absolutely continuous with respect to the Lebesgue measure, it follows from  Thaler (1980) that:
\begin{itemize}
\item[-] if $\gamma \in ]0,1[$,
there exists a
unique absolutely continuous $T_\gamma$-invariant probability measure $\nu_\gamma$ (or simply $\nu$)
on $[0, 1]$;
\item[-] if $\gamma\geq 1$,
there is no absolutely continuous invariant probability measure.
\end{itemize}

\begin{figure}[htb]
\centerline{
\includegraphics[width=7.2cm, trim=33mm 90mm 40mm 80mm]{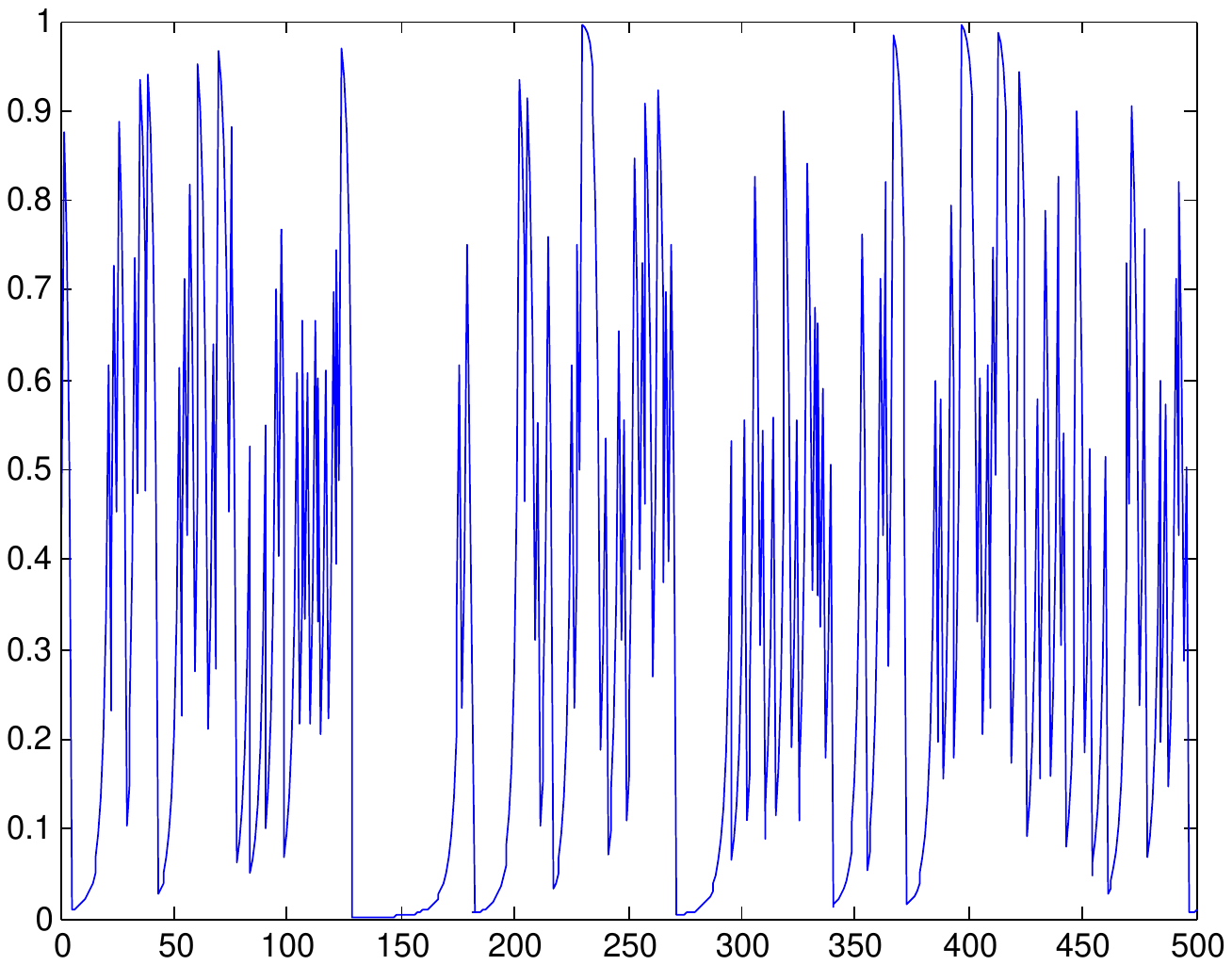}
\includegraphics[width=7.2cm, trim=40mm 90mm 33mm 80mm]{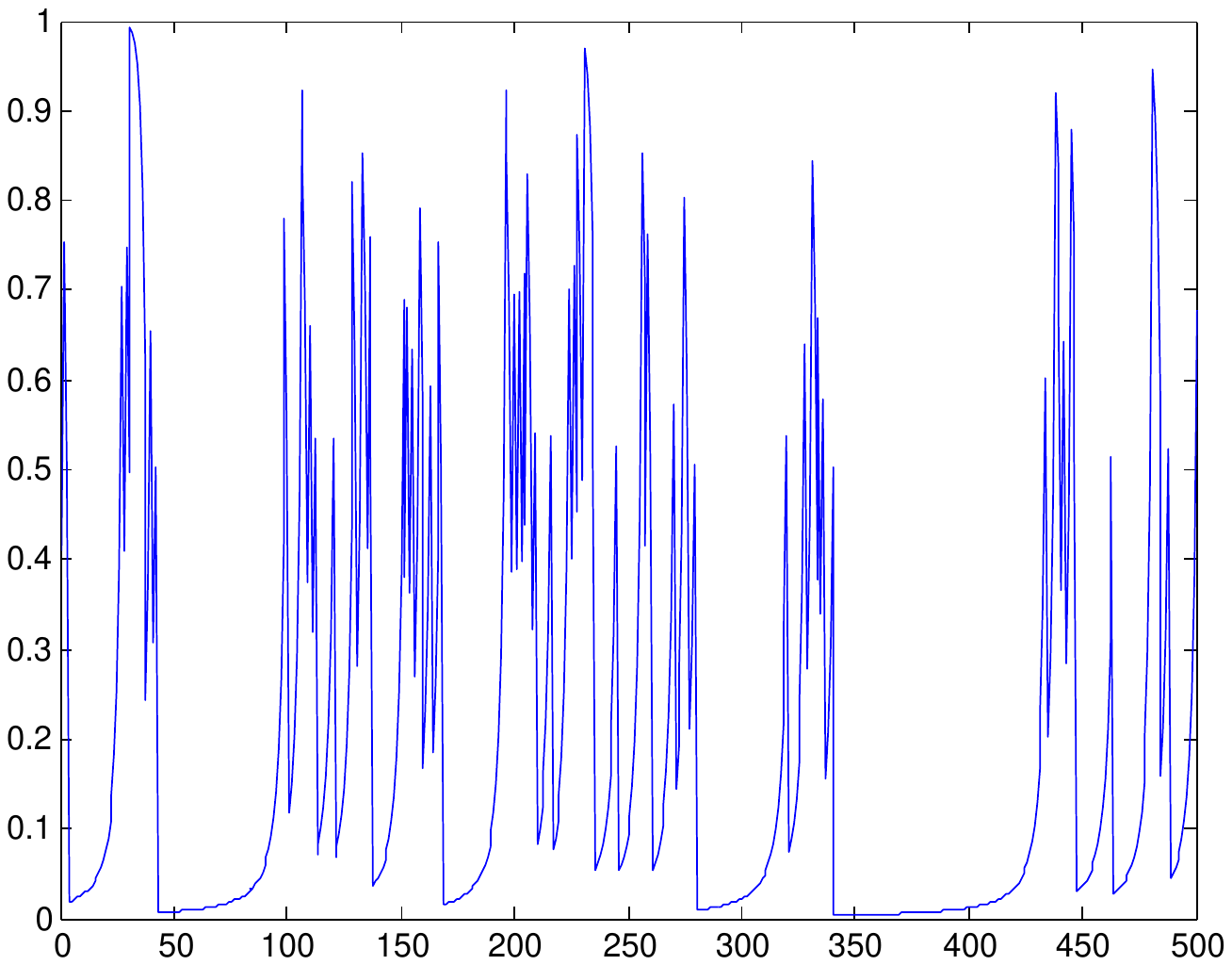}
}
\caption{Time series of $500$ iterations of the  intermittent map $T_{0.5}$ (left) and $T_{0.9}$ (right).}
\label{fig:pm_ts}
\end{figure}
For $x$ near the neutral fixed point $0$,  the sequence $T^k(x), k\geq 0$, spends
a lot of time around  $0$, since $T(x)=x(1+(2x)^\gamma) \approx x$ for $x\approx 0$; see Figure~2 for the time series of $500$ iterations of $T_{0.5}$ (left) and $T_{0.9}$ (right). Note that the length of the periods spent in the neighborhood of zero increases as $\gamma$ gets larger.
Hence the
density $h_\gamma$
(or simply $h$) of the invariant distribution $\nu_\gamma$ explodes in
the neighborood of $0$.
Even though no explicit formula is known for $h$, we can give a precise description thanks to the works by Thaler (1980) and Liverani, Saussol and Vaienti (1999). From the first paper, we infer  that  the function  $x \mapsto x^\gamma h(x)$ is bounded from above and below.
From the second paper, we know that $h$ is non increasing with $h(1)>0$, and that it is Lipshitz on any interval $[a,1]$ with $a>0$.
Since $h$ is strictly positive
on $[0,1]$,
the probability $\nu$ is equivalent to the Lebesgue measure on $[0,1]$.

For $x\in [\epsilon,1]$, the map $T$ is expanding, i.e. $|T^\prime(x)| \geq \alpha$, for some $\alpha=\alpha(\epsilon)>1$, and thus the sequence $T^k(x)$ exhibits random behavior.  As
$T(1/2+\epsilon)=2\epsilon$, the sequence occasionally returns to a neighborhood of $x=0$. Then, as
$T(x)\approx x$ for $x\approx 0$, the sequence $T^k(x)$ hovers around the neighborhood of zero for a long time, which may explain the long range dependence of the process $T^k$, $k\geq 1$.
Eventually, since $T(x)=x+2^\gamma x^{\gamma +1}>x$, and
$T^\prime(x) =1+2^\gamma (\gamma+1) \gamma x^{\gamma} >1$, the process escapes the neighborhood of zero and reenters the zone of chaotic behaviour.

From now on, we shall use the notation
$$
\nu(f)= \int f(x) \nu(dx) \, ,
$$
which is valid for any $f\in {\mathbb L}^1([0,1], \nu)$.
For $\gamma \in ]0,1[$, we view $T^k$ as a random variable from the probability space $([0,1], \nu)$ to $[0,1]$.
The fact that $\nu$ is invariant by  $T$  implies that $\nu(f)=\nu(f \circ T^k)$ for any $f\in {\mathbb L}^1([0,1], \nu)$, and more generally, it also implies
that the process $(T^k)_{k \geq 0}$ is strictly stationary.
Let us  briefly recall some known results about the iterates of $T$:
\begin{enumerate}
\item {\em Decay of correlations.} For any bounded function $f$ and any H\"older functions $g$, Young (1999) proved the following decay  of the covariances
\begin{equation}\label{e:rate}
  \nu\Big((g-\nu(g))\cdot (f-\nu(f))\circ T^n  \Big)=
  \nu \Big((g-\nu(g))\cdot f \circ T^n\Big) = O(n^{(\gamma-1)/\gamma}),
\end{equation}
as $n\rightarrow \infty$.
Some lower bounds for the covariance in \eref{e:rate} can be found in the paper by Sarig (2002), proving that the rate $n^{(\gamma-1)/\gamma}$
is optimal. Dedecker, Gou\"ezel and Merlev\`ede (2010) have shown that \eref{e:rate}  remains true if $g$ is
any bounded variation function.
\item {\em Behaviour of Birkhoff sums.} Liverani, Saussol and Vaienti (1999) have proved that the map $T$ is mixing in the ergodic
theoretic sense. Let then
$$
s_n(f)= \sum_{k=1}^n f \circ T^k \, .
$$
For any $f\in {\mathbb L}^1([0,1], \nu)$, it follows
from Birkhoff's ergodic theorem that $n^{-1}s_n(f)$ converges to $\nu(f)$ almost everywhere.

Concerning the convergence in distribution of the sequence $s_n(f)-n \nu(f)$ (suitably normalized) on the probability space $([0,1], \nu)$, we must distinguish three cases. If $\gamma \in ]0,1/2[$
and $f$ is any H\"older function,  Young (1999) proved that the sequence  $n^{-1/2}(s_n(f)-n\nu(f))$ converges in distribution to a normal law.
Next,  Gou\"ezel (2004) has given a complete picture of the convergence in distribution of $s_n(f)-n\nu(f)$ when $\gamma \in [1/2, 1[$.
More precisely, if  $\gamma =1/2$  and $f$ is any  H\"older function, he proved that the sequence  $(n\log(n))^{-1/2}(s_n(f)-n\nu(f))$ converges in distribution to a normal law.
If $\gamma \in ]1/2,1[$ and $f$ is any H\"older function, he proved that $n^{-\gamma}(s_n(f)-n \nu(f))$ converges in distribution to a stable law of index $1/\gamma$.

\end{enumerate}

The power decay of the covariance \eref{e:rate} suggests that there may be long-range dependence
for some values of $\gamma$. A finite variance stationary process with covariance $r(n)$ is said to be {\it short-range dependent}
if $\sum_{n=0}^\infty |r(n)|< \infty$ and {\it long-range dependent} if
$\sum_{n=0}^\infty |r(n)|=\infty$. In view of the optimality of \eref{e:rate}, we see that the process
$(T^k)_{k \geq 0}$ is short-range dependent if $\gamma \in ]0, 1/2[$ and long-range dependent if $\gamma \in ]1/2, 1[$. The case $\gamma=1/2$ is a boundary case.
Moreover, the asymptotic behaviour of the normalized sums $s_n(f)$ is normal in the short-range dependent case (including $\gamma=1/2$), and is stable in the long-range dependent case.

Our aim is  to
study the limit in distribution  of the empirical process
\begin{equation}\label{e:G}
G_n(t)=\frac 1 n \sum_{k=1}^n
\Big({\mathbf 1}_{T^k \leq t}-F(t)\Big), \quad t \in [0, 1],
\end{equation}
with $F(t)= \nu([0,t])$,
in the case where $\gamma \in [1/2, 1[$.

Let us introduce another stationary process with the same law.
Let first $K$ be the
 Perron-Frobenius operator of $T$ with respect to $\nu$, defined as follows:
for any  functions $f, g$ in ${\mathbb L}^2([0,1], \nu)$
\begin{equation}\label{Perron}
\int  f (T(x)) g(x) \nu(dx)= \int f(x) (K g) (x)  \nu(dx) \, .
\end{equation}
The relation \eref{Perron} states  that $K$ is  the adjoint operator of the isometry $U: f \mapsto f\circ T$
acting on ${\mathbb L}^2([0,1], \nu)$. It is easy to see that the operator $K$ is a transition kernel.\footnote{Indeed, by stationarity, the relation \eref{Perron}
can be written as follows:
for any  functions $f, g$ in ${\mathbb L}^2([0,1], \nu)$,
$$
\nu((f \circ T)\cdot g)=\nu((f \circ T) \cdot (K g) \circ T)\, .
$$
On the probability space $([0,1], \nu)$, this means precisely that $(Kg)\circ T={\mathbb E}(g|T)$.
Hence
\[
(Kg)(x)={\mathbb E}(g|T=x),
\]
so that $K$ is a transition kernel.}
Let now $(X_i)_{i \in {\mathbb Z}}$ be a stationary Markov chain with invariant measure $\nu$
and transition kernel $K$. It is well known (see for instance Lemma XI.3 in
Hennion and Herv\'{e} (2001)) that on the probability space $([0, 1], \nu)$, the
random vector $(T, T^2, \ldots , T^n)$ is distributed as
$(X_n,X_{n-1}, \ldots, X_1)$.\footnote{For instance, by the Perron-Frobenius relation (\ref{Perron}) and
stationarity $${\mathbb E}(f(X_1)g(X_2))=\nu(f\cdot Kg)=\nu((f \circ T)\cdot g)=\nu(f(T^2) g(T)).$$
By setting $f(x)=\exp(itx)$ and $g(y)=\exp(iuy)$, we obtain that $(X_1,X_2)$ is distributed as $(T^2,T)$. }
Hence, the process
\begin{equation}\label{e:L}
L_n(t)=\frac 1 n \sum_{k=1}^n
\Big({\mathbf 1}_{X_k \leq t}-F(t)\Big), \quad t \in [0,1]
\end{equation} has the same distribution as
$\{G_n(t), t \in [0,1]\}$.

In the short-range dependent case $\gamma \in ]0, 1/2[$, Dedecker (2010)  proved that,
on the probability
space $([0, 1], \nu_\gamma)$ the process $\{\sqrt{n} G_n(t), t \in [0,1]\}$
converges in distribution in the space $D([0, 1])$ of cadlag functions equipped with the uniform metric to a centered Gaussian process
$G$, whose sample paths are almost surely uniformly continuous. Moreover the
covariance function of $G$ is given by
\begin{equation}\label{covemp2}
\mathrm{Cov}(G(s), G(t))=\nu (f_t^{(0)} \cdot f^{(0)}_s) + \sum_{k>0} \nu (f_t^{(0)} \cdot f_s^{(0)} \circ T^k) +
\sum_{k>0} \nu (f_s^{(0)} \cdot f_t^{(0)} \circ T^k) \, ,
\end{equation}
where the function $f^{(0)}_t$ is defined by
$$
f^{(0)}_t(x)={\bf 1}_{x \leq t} -
\nu([0, t])\, .
$$
For $s=t$, the series  \eqref{covemp2} is the asymptotic variance of $\sqrt{n} G_n(t)$.
This variance  has the same structure as  the asymptotic variance of the
normalized partial
sums of a stationary sequence $(Y_i)_{i \geq 0}$ in the case where the covariance
series converge, that is
$$
{\mathrm{Var}}(Y_0)+ 2\sum_{k=1}^\infty {\mathrm{Cov}}(Y_0, Y_k)\, .
$$
Observe that in this case, the limit process $\{G(t), t \in [0,1]\}$ is not degenerate.

In the long range-dependent case  $\gamma \in [1/2, 1[$,  the series in (\ref{covemp2}) may not
converge.
The
long-range dependent case   has been studied  by Dehling and Taqqu (1989) for the empirical process
of a stationary Gaussian sequence. In that paper, the authors show that the empirical process, suitably normalized,
converges in distribution in $D({\mathbb R})$ to a degenerate Gaussian process.
Following the approach of Dehling and Taqqu, Surgailis (2002)  proved that the empirical process
of a linear process whose innovations belong to the domain of normal attraction of a stable distribution,
converges  in distribution in $D({\mathbb R})$  to a degenerate stable process.
In the two papers cited above, the main idea is to approximate the empirical process by a sum
of independent random variables whose asymptotic distribution is easy to derive. Such an approximation
is not available in our context, and we shall use a completely different approach.
For $\gamma \in [1/2, 1[$, we shall obtain the same limit behavior as in Surgailis (2002), but
in a space whose topology is much coarser than that of $D([0,1])$.


We shall investigate here the behavior of the empirical
process
$$\{G_n(t), t \in [0,1]\} \quad \text{for  $\gamma \in [1/2,1[$,}$$
in the
Hilbert space
 $$ H={\mathbb L}^2\big([0,1], dt\big)
 $$
 with norm $\|\cdot\|_H$.   The reason is that we can use very precise
deviation inequalities for $H$-valued random variables to prove the
tightness of the empirical process in $H$. An interesting question is
whether our results remain true in $D([0, 1])$, as in
Dehling and Taqqu
(1989) and Surgailis (2002).

Note that the empirical process, viewed as element of the Hilbert space $H$, is a centered and normalized
partial sum of the
random variables $\xi_i$, defined by
\[
 \xi_i(t)=1_{[0,t]}(X_i).
\]
If the underlying random variables $(X_i)_{i\geq 1}$ are independent and indentically distributed (i.i.d.) or mixing, weak convergence of the empirical process is an immediate corollary of an appropriate central limit theorem for $H$-valued random variables. Such CLTs have been established in the i.i.d. case by Mourier (1953), for $\phi$-mixing processes by Kuelbs and Philipp (1980) and for strongly mixing processes by Dehling (1983).
The problem in our situation is that the process $(X_k)_{k\geq 1}$ is not mixing, so that none of the known CLTs for Hilbert space valued random variables can be applied. In this paper, we will present a proof that is
tailor-made for the  empirical process of intermittent maps.

Concerning the weak convergence of the empirical process
$\{G_n(t), t \in [0,1]\}$ with $\gamma \in [1/2, 1[$
in the space  $H$,
we shall prove the following theorem.

\begin{thm}\label{mainth}
On the  probability space $([0,1], \nu)$, the following results hold:
\begin{enumerate}
\item If $\gamma=1/2$, then  the process
$$\Big\{\frac{\sqrt{n}}{\sqrt {\log (n)}} G_n(t), t \in [0,1]\Big\}$$ converges in distribution in $H$ to a degenerate Gaussian process
$\{  g(t) Z, t \in [0,1]\}$ where $Z$ is a standard normal
and $g(t)=\sqrt {h(1/2)} (1-F(t))$.
\item If $\gamma \in ]1/2, 1[$, then  the process
$$\{n^{1-\gamma}G_n(t), t \in [0,1]\}$$ converges in distribution in $H$ to a degenerate stable  process
$\{  g(t) Z, t \in [0,1]\}$ where
 $g(t)= C_\gamma (h(1/2)) ^\gamma
(1-F(t))$
with
$$
C_\gamma= \frac{1}{4^\gamma \gamma}
\Big ( \Gamma(1-1/\gamma)\cos\Big(\frac{\pi}{2\gamma}\Big)\Big)^\gamma \, ,
$$
and $Z$ is an
$1/\gamma$-stable random variable totally skewed to the right, that is with characteristic function
\begin{equation}\label{char}
{\mathbb E}(\exp(itZ))=\exp\Big(-|t|^{1/\gamma}\Big(1-\mathrm{sign}(t) \tan
\Big(\frac{\pi}{2\gamma}\Big)\Big)\Big)\, .
\end{equation}
\end{enumerate}
\end{thm}

\begin{rem}
Recall that, on the probability space $([0,1], \nu)$ the process
$\{G_n(t), t \in [0,1]\}$ is distributed as
$\{L_n(t), t \in [0,1]\}$, defined in \eref{e:L}.
Hence Theorem \ref{mainth} is also valid for the
empirical process $\{L_n(t), t \in [0,1]\}$.
\end{rem}

\begin{rem}
Recall also that the distribution $S_\alpha(\sigma, \beta, \mu)$ of a stable
random variable is characterized by the stability parameter $\alpha \in ]0, 2]$,
the scale parameter $\sigma>0$, the skewness parameter $\beta \in [-1,1]$
and the shift parameter $\mu \in {\mathbb R}$ (see e.g. Samorodnitsky and Taqqu (1994)). When $ \gamma \in ]1/2, 1[$, the random
variable $Z$ in \eref{char}  has a stability parameter $\alpha=1/\gamma \in ]1, 2[$ and hence has infinite variance and finite mean. Moreover $\sigma=1, \beta=1$ and $\mu=0$. Since
$\beta=1$ it is said to be ``totally skewed to the right".
\end{rem}
\begin{rem}
The limits in the short-range and long-range dependent case are quite different.
As noted after the relation \eref{e:L}, in the short-range dependent case, the limit is a
non-degenerate Gaussian process, whereas in the long-range dependent case considered
in Theorem \ref{mainth} the limit is a degenerate process
$\{g(t) Z, t \in [0,1]\}$, where $Z$ can be Gaussian or not depending on
the value of $\gamma$.
\end{rem}

Concerning the almost sure behavior of
$\|G_n\|_H$ , we shall prove the following theorem.
\begin{thm}\label{as}
The following results hold:
\begin{enumerate}
\item Let $\gamma=1/2$.
Let $a_n$ be any  sequence of numbers  such that $a_n\geq a$
for some $a>0$, and
$\sum_{n>0} n^{-1} a_n^{-2}< \infty$.
Then for any $\varepsilon>0$,
\begin{equation}\label{dev1}
\sum_{n=1}^\infty \frac{1}{n}
{\mathbb \nu}\Big ( \max_{1 \leq k \leq n} \frac{k\|G_k\|_H}
{\sqrt{n \log (n+1)} a_n} \geq \varepsilon
\Big ) < \infty \, .
\end{equation}
Assume moreover that $a_n$ is non decreasing  and such that
$a_n \leq c a_{n/2}$ for some $c\geq 1$. Then
\begin{equation} \label{e:asG1}
  \frac{\sqrt{n}}
{ a_n \sqrt{\log(n)}}  \|G_n\|_H \quad \text{converges almost everywhere to 0}\,  .
\end{equation}
\item Let $\gamma \in ]1/2,1[$.
Let $a_n$ be any  sequence of positive numbers such that
$\sum_{n>0} n^{-1} a_n^{-1/\gamma}< \infty$.
Then for any $\varepsilon>0$,
\begin{equation}\label{dev2}
\sum_{n=1}^\infty \frac{1}{n}
{\mathbb \nu}\Big ( \max_{1 \leq k \leq n} \frac{k\|G_k\|_H}
{n^\gamma a_n} \geq \varepsilon
\Big ) < \infty \, .
\end{equation}
Assume moreover that $a_n$ is non decreasing  and such that
$a_n \leq c a_{n/2}$ for some $c\geq 1$. Then
\begin{equation} \label{e:asG2}
  \frac{n^{1-\gamma}}
{ a_n}  \|G_n\|_H \quad \text{converges almost everywhere to 0}\,  .
\end{equation}
\end{enumerate}
\end{thm}
\begin{rem}
The corresponding almost sure result is true also for
$\|L_n\|_H $, see the proof of Theorem \ref{as}.
\end{rem}
\begin{rem}
For instance, all the conditions on $a_n$ are satisfied if $a_n= (\log(n+1))^\delta$ for $\delta>\gamma$.
For $\gamma \in ]1/2, 1[$  this is in  accordance with the i.i.d. situation,
which we now describe.

 Let $(X_i)_{i \geq 1}$ be a sequence of
i.i.d. centered random variables such that  $n^{-\gamma} (X_1+
\cdots + X_n)$ converges in distribution to the $1/\gamma$-stable distribution
with characteristic function \eref{char}. It is well known (see for instance
Feller (1966), page 547) that this  implies that $x^{1/\gamma}
{\mathbb P}(X_1< -x)\rightarrow 0$ and $x^{1/\gamma} {\mathbb
P}(X_1>x)\rightarrow c>0$ as $x \rightarrow \infty$. For any
nondecreasing sequence $(b_n)_{n \geq 1}$ of positive numbers,
either $(X_1+ \cdots + X_n)/b_n $ converges to zero almost
surely or $\limsup_{n \rightarrow \infty} |X_1+ \cdots +
X_n|/b_n=\infty$ almost surely, according as $\sum_{n=1}^\infty
{\mathbb P}(|X_1|>b_n)< \infty$ or $\sum_{n=1}^\infty {\mathbb
P}(|X_1|>b_n)= \infty$ -- this follows from the proof of
Theorem 3 in Heyde (1969). If one takes $b_n=n^{\gamma}(\ln(n+1))^\delta$
we obtain the constraint $\delta>\gamma$ for the almost sure
convergence of $n^{-\gamma}(\ln(n+1))^{-\delta} (X_1+ \cdots +X_n)$ to
zero. This is exactly the same constraint as in our dynamical
situation.

This situation is similar to the one described in Theorem 1.7 of Dedecker,
Gou\"ezel ans Merlev\`ede (2010).
Note that there is a mistake in Theorem 1.7 of this paper, in the case where
$p=1/2$ (weak moment of order 2): the exponent of the logarithm in (1.8) should satisfy $b>1$ instead of $b>1/2$.
\end{rem}
\begin{rem}
In the short-range dependent case $\gamma \in ]0, 1/2[$, Dedecker, Merlev\`ede and Rio (2013) have proved
a strong approximation result for the empirical process $\{G_n(t), t \in [0,1]\}$.
As a consequence, it follows that: almost everywhere, the sequence
$$
\Big \{\frac{\sqrt {n}}{\sqrt {2\log \log (n)}} G_n(t), t \in [0,1] \Big \}
$$
is relatively compact for the supremum norm, and
the set of limit points is the unit ball of the reproducing kernel
Hilbert space associated with the covariance function  \eqref{covemp2} of the limit Gaussian proces $G$.  In particular, it follows that
$$
\limsup_{n \rightarrow \infty} \frac{\sqrt {n}}{\sqrt {2\log \log (n)}} \|G_n\|_H=
\sigma  \quad \text{almost everywhere, where
$\sigma= \sup_{\|g\|_H \leq 1} \sqrt{{\mathbb E}(<g, G>^2)}\, .
$}
$$
\end{rem}

\medskip

To conclude this section, let us note that the conclusions of Theorems
 \ref{mainth} and \ref{as} also hold when $T^k$ is replaced by $g(T^k)$ where
$g$ is a monotonic and H\"older continuous function
from $[0,1]$ to ${\mathbb R}$. This fact  will be proved in
Theorem \ref{prop:gene} of Section \ref{Sec:new}, and used in Section \ref{Sec:CVM}.

\section{Application to the weak convergence of the  Wasserstein distance $W_1$.}
\setcounter{equation}{0}
Let $\mu_n$ be the empirical measure of the iterates of $T$, that is
$$
\mu_n= \sum_{k=1}^n \delta_{T^k} \, .
$$
Consider the  Wasserstein distance between the empirical measure
$\mu_n$ and $\nu$:
$$
W_1(\mu_n, \nu)=\inf_{ \pi \in M(\mu_n, \nu) } \int |x-y| \pi (dx, dy)\, ,
$$
where $M(\mu_n, \nu)$ is the set of probability measures
on $[0,1]^2$ with marginals $\mu_n$ and $\nu$. Since $\mu_n$
and $\nu$ are probability measures
on the real line, it is well known that (see for instance  Fr\'echet (1957))
\begin{equation}\label{W1}
W_1(\mu_n, \nu)= \int_0^1 |G_n(t)|dt \, ,
\end{equation}
where $G_n$ is defined in \eref{e:G}.
Since the functional $\psi(f)=\int_0^1 |f(t)| dt$ is  continuous
on the Hilbert space $H$, we can apply Theorem \ref{mainth}.
Since
$$
\int_0^1 (1-F(t))dt = \int_0^1 x \nu(dx)=  \int_0^1 x h(x) dx \, ,
$$
we obtain
\begin{cor}
On the  probability space $([0,1], \nu)$, the following results hold:
\begin{enumerate}
\item If $\gamma=1/2$,
$$
\frac{\sqrt{n}}{\sqrt{\log(n)}} W_1(\mu_n, \nu) \quad
\text{converges in distribution to} \quad
\sqrt{h(1/2)} |Z| \int_0^1 x h(x) dx \, .
$$
where $Z$ is a standard normal.
\item If $\gamma\in ]1/2,1[$,
$$
n^{1-\gamma} W_1(\mu_n, \nu) \quad
\text{converges in distribution to} \quad
C_\gamma (h(1/2))^\gamma|Z| \int_0^1 x h(x) dx \, .
$$
where $Z$ is an
$1/\gamma$-stable random variable with characteristic function
(\ref{char}).
\end{enumerate}
\end{cor}
Since by \eqref{W1}, $W_1(\mu_n, \nu) \leq \|G_n\|_H$, we can apply
Theorem \ref{as}.
\begin{cor}
The following results hold:
\begin{enumerate}
\item Let $\gamma=1/2$
and let  $a_n$ be as in Item 1 of Theorem \ref{as}. Then
$$
\frac{\sqrt n} {a_n \sqrt{\log(n)}} W_1(\mu_n, \nu) \quad
\text{converges  to zero almost everywhere.}
$$
\item
Let $\gamma \in ]1/2, 1[$, and let  $a_n$ be as in Item 2 of Theorem \ref{as}. Then
$$
\frac{n^{1-\gamma}} {a_n} W_1(\mu_n, \nu) \quad
\text{converges  to zero almost everywhere.}
$$
\end{enumerate}
\end{cor}

\section{Application to the Cram\'er-von-Mises statistic}\label{Sec:CVM}
Recall that $F$ is the cumulative distribution function (cdf)
of the absolutely continuous invariant measure $\nu$.
In order to test whether the cdf of $\nu$ is equal to the cdf $G$, we use the test statistic
$$
 \Psi_n=\int_0^1 (F_n(t)-G(t))^2 dG(t),
\quad \text{where} \
F_n(t)= \frac{1}{n} \sum_{i=1}^n {\bf 1}_{T^i \leq t} \, .
$$
Note that, by Birkoff's
ergodic  Theorem
combined with the Glivenko-Cantelli arguments,
$$
   \sup_{t \in [0,1]} |F_n(t)-F(t)| \quad \text{converges almost everywhere to 0.}
$$
It follows that
$\Psi_n$ converges almost everywhere to $0$ if $G=F$ and to
a strictly positive number
if $G$ is absolutely continuous and $G\neq F$.

We can apply Proposition \ref{prop:gene} of the Appendix (which is a generalization of
Theorem \ref{mainth}) to derive the asymptotic distribution of $\Psi_n$
on the probability space $([0,1], \nu)$,
under the null hypothesis
$
H: F=G
$.
\begin{cor}
\label{cor:cvm-stat}
On the probability space $([0,1], \nu)$, under the null
hypothesis $
H: F=G$, the following results hold:
\begin{enumerate}
\item
If  $\gamma=1/2$,
\[
  \frac{n}{\log n} \Psi_n \quad
\text{converges in distribution to} \quad \frac{1}{3} h(1/2)  Z^2,
\]
where $Z$ is a standard normal random variable.
\item
If $\gamma\in ]1/2,1[$,
\[
  n^{2-2\gamma} \Psi_n \quad
\text{converges in distribution to} \quad \frac{C_\gamma^2 ( h(1/2)  )^{2\gamma}}{3}
   Z^2,
\]
where $Z$ is a $1/\gamma$-stable random variable with characteristic function (\ref{char}).
\end{enumerate}
\end{cor}
\noindent{\em Proof.} Assume that $F=G$, and note that $F$ is continuous and strictly increasing.
By a change of variables, we obtain
\[
 \Psi_n=\int_0^1 (G_n(t))^2 \nu(dt) = \int_0^1 (G_n(F^{-1}(t)))^2 dt.
\]
Since  $F^{-1}(F(t))=t$,
\[
G_n(F^{-1}(t)) = \frac{1}{n}\sum_{i=1}^n \Big( {\bf 1}_{T^i\leq F^{-1}(t)}-t\Big)=\frac{1}{n}\sum_{i=1}^n \Big( {\bf 1}_{F(T^i)\leq t} -t \Big),
\]
i.e. $G_n(F^{-1}(t))$ is the empirical process  of the sequence
$(F(T), F(T^2), \ldots, F(T^n))$.

In addition of being increasing, the function $F$ is also H\"older continuous. Indeed,
let $\delta \in ]0, 1-\gamma[$ and $p=1/(1-\delta)$. By H\"older's inequality,
for any $x,y \in [0,1]$ with $x\leq y$,
$$
  F(y)-F(x) \leq \Big( \int_x^y h(x)^p dx\Big )^{1/p} (y-x)^\delta \, .
$$
Note that $\int_0^1 h(x)^p dx$ is finite, because $x^\gamma h(x)$ is bounded from above  and $p\gamma<1$. It follows that there exists a positive constant
$K$ such that
$$
F(y)-F(x) \leq K(y-x)^\delta \, ,
$$
showing that $F$ is H\"older continuous of index $\delta$.

Hence one can apply Theorem \ref{prop:gene} of Section \ref{Sec:new}: on the
probability space $([0,1], \nu)$,
\begin{enumerate}
\item If $\gamma=1/2$, then  the process
$$\Big\{\frac{\sqrt{n}}{\sqrt {\log (n)}} G_n(F^{-1}(t)), t \in [0,1]\Big\}$$ converges in distribution in $H$ to a degenerate Gaussian process
$\{  g(t) Z, t \in [0,1]\}$ where $Z$ is a standard normal
and $g(t)=\sqrt {h(1/2)} (1-t)$.
\item If $\gamma \in ]1/2, 1[$, then  the process
$$\{n^{1-\gamma}G_n(F^{-1}(t)), t \in [0,1]\}$$ converges in distribution in $H$ to a degenerate stable  process
$\{  g(t) Z, t \in [0,1]\}$ where
 $g(t)= C_\gamma (h(1/2)) ^\gamma
(1-t)$
and $Z$ is a $1/\gamma$-stable random variable with characteristic function (\ref{char}).
\end{enumerate}
 Since the map $f\mapsto \|f\|_H^2=\int_0^1 f^2(t)dt$ is
continuous on $H$, we obtain in case $\gamma=1/2$ that
\[
  \frac{n}{\log n}\Psi_n=\left\| \frac{\sqrt{n}}{\sqrt{\log n}}  G_n \circ F^{-1} \right\|_H^2
\text{converges in distribution to} \
 h(1/2) Z^2\int_0^1 (1-t)^2 dt =\frac{1}{3} h(1/2) Z^2.
\]
The case $1/2<\gamma < 1$ follows in the same way.

\section{Proof of Theorems  \ref{mainth} and \ref{as}}
In this section, $C$ is a positive constant which may vary from
line to line.
\setcounter{equation}{0}
\subsection{Some general facts}\label{Sec:gen}

Let
\begin{equation}\label{e:Yt}
Y_i(t)= {\mathbf 1}_{X_i \leq t}-F(t) \,  ,
\end{equation} and
$S_n=\sum_{i=1}^n Y_i$.  With these notations, by \eref{e:L},
\begin{equation}\label{e:Sn}
S_n(t)= \sum_{i=1}^n ({\mathbf 1}_{X_i \leq t}-F(t))= n L_n(t)\, .
\end{equation}
Let also  $V_i(t)= {\mathbf 1}_{T^i \leq t}-F(t)$, and
$\Sigma_n=\sum_{i=1}^n V_i$.  With these notations, by \eref{e:G},
\begin{equation}\label{e:Sigman}
\Sigma_n(t)= \sum_{i=1}^n ({\mathbf 1}_{T^i \leq t}-F(t))= n G_n(t)\, .
\end{equation}

Recall that, on the probability space $([0,1], \nu)$ the
sequence $(V_1, V_2, \ldots, V_n)$ is distributed as the sequence
$(Y_n, Y_{n-1}, \ldots, Y_1)$. It follows that $L_n$ is
distributed as $G_n$, and it is equivalent to prove
Theorem \ref{mainth} for $L_n$ or for  $G_n$.

\begin{enumerate}
\item
Let us first prove the following inequality:
 for any $x
\geq 0$,
\begin{equation}
\label{equ2law}
\nu \Big( \max_{1 \leq k \leq n} \|\Sigma_k\|_H \geq x \Big) \leq
{\mathbb P}\Big(2\max_{1 \leq k \leq n} \|S_k\|_H \geq x \Big) \,  .
\end{equation}
Indeed,
\begin{equation}
  \label{equ1law}
  \max_{1 \leq k \leq
  n} \Big \| \sum_{i=1}^{k} V_i \Big \|_H \quad \text{is distributed as} \quad  \max_{1 \leq k \leq
  n} \Big \| \sum_{i=k}^n X_i \Big \|_H \, .
\end{equation}  Notice now that for any $k \in \{1, \ldots, n\}$,
\[
\sum_{i=k}^n X_i = \sum_{i=1}^n X_i - \sum_{i=1}^{k-1} X_i \, .
\]
Consequently
\[
\max_{1 \leq k \leq
  n} \Big\|\sum_{i=k}^n X_i \Big\|_H \leq  \max_{1 \leq k \leq
  n-1} \Big \|\sum_{i=1}^{k}X_i\Big\|_H + \Big \|\sum_{i=1}^{n}X_i\Big\|_H
  \leq 2\max_{1 \leq k \leq
  n} \Big \|\sum_{i=1}^{k}X_i\Big\|_H
\]
which together with \eqref{equ1law} entails \eqref{equ2law}.
\item
Let $H$ be a separable Hilbert space with inner product $\langle x, y\rangle$ and norm $\|x\|_H^2=\langle x,x\rangle$, and let $(u_i)_{i\geq 1}$ be a complete orthonormal system in $H$. Thus, any vector $x\in H$ can be expanded into a series
\[
  x=\sum_{i=1}^\infty \langle x,u_i\rangle u_i.
\]
For any integer $m$, we define the finite-dimensional projection $P_m:H\rightarrow H$ by
\[
 P_m(x) =\sum_{i=1}^m \langle x,u_i\rangle u_i.
\]
Let now $(Z_n)_{n\geq 1}$ be a sequence of $H$-valued random variables, and let $Z$ be another $H$-valued random variable, satisfying
\\[.5mm]
(i) For all integers $m$, as $n\rightarrow \infty$,
\begin{equation}
 \text{$P_m(Z_n)$ converges in distribution in $H$ to $P_m(Z)$}
\label{eq:fidi-1}
\end{equation}
(ii) For any $\epsilon >0$,
\begin{equation}
 \lim_{m\rightarrow \infty} \limsup_{n\rightarrow \infty} P(\|Z_n-P_m(Z_n)\|_H\geq \epsilon)=0.
\label{eq:fidi-approx}
\end{equation}
Then, by Theorem~4.2 of Billingsley (1968), $Z_n$ converges in distribution to $Z$, as $n\rightarrow \infty$.

Observe that weak convergence in \eqref{eq:fidi-1} is really weak convergence in a finite dimensional Euclidean space, since the map $a=(a_1,\ldots,a_m) \mapsto \sum_{i=1}^m a_i u_i$ defines an isometry between ${\mathbb R}^m$ and the subspace $P_m(H)\subset H$. Thus \eqref{eq:fidi-1} holds if and only if
\[
  \left(\langle Z_n, u_1 \rangle,\ldots, \langle Z_n, u_m \rangle  \right)
  \  \text{converges in distribution to}  \
\left(\langle Z, u_1 \rangle,\ldots, \langle Z, u_m \rangle  \right),
\]
as $n\rightarrow \infty$. Hence, using the Cram\'er-Wold device, \eqref{eq:fidi-1} is equivalent to
\begin{equation}
\langle Z_n,\sum_{i=1}^m \lambda_i u_i \rangle = \sum_{i=1}^m \lambda_i \langle Z_n,u_i \rangle
\  \text{converges in distribution to}  \
\sum_{i=1}^m \lambda_i \langle Z ,u_i \rangle = \langle Z, \sum_{i=1}^m \lambda_i u_i \rangle,
\end{equation}
as $n\rightarrow \infty$, for all $\lambda_1, \ldots, \lambda_m\in {\mathbb R}$.

\item
On the Hilbert space $H={\mathbb L}^2([0,1], dt)$, consider the usual orthonormal Fourier basis
$$({\bf 1}, (\sqrt 2 \cos(2\pi k \cdot), \sqrt 2 \sin(2\pi k\cdot))_{k \in {\mathbb N}^*})\, ,$$ and let $P_m$ be the projection operator on the space
spanned by $$({\bf 1}, (
\sqrt 2 \cos(2\pi k \cdot), \sqrt 2 \sin(2\pi k\cdot))_{k \in
\{1, \ldots, m\}})\, .$$
Let $\|\cdot \|_\infty$ denotes the essential supremum norm, namely
$$\|X\|_\infty= \inf \{ M>0: {\mathbb P}(|X|>M)=0\}\, .$$
Thus while $\| Y_0-P_m(Y_0) \|_H$ is random,
\begin{equation}\label{e:Mm}
M(m)=\|\| Y_0-P_m(Y_0) \|_H \|_\infty
\end{equation}
is not random. In fact, we have
\begin{equation}\label{e:supnorm}
M(m) = O(m^{-1/2})\, .
\end{equation}
Indeed, the Fourier coefficients
 $$a_k(x)= \sqrt 2 \int_x^1 \cos (2\pi kt) dt \quad \text{and} \quad b_k(x)=
 \sqrt 2 \int_x^1 \sin (2\pi kt) dt $$
 of the functions $f_x(t)= {\bf 1}_{t\geq x}$ are such that
 $$
 a_k(x)+ i b_k(x)=\sqrt{2} \int_x^1\exp (i2\pi kt) dt= \frac{1}{i \sqrt{2} \pi k}(1-\exp (i2\pi kx))\, ,
 $$
 in such a way that
 $a_k^2(x)+b_k^2(x)\leq 2/(\pi k)^2$.
 Hence,
\begin{equation}\label{e:bproj}
\| Y_0-P_m(Y_0) \|^2_H= \sum_{k=m+1}^\infty (a_k^2(X_0)
+ b_k^2(X_0)) \leq \sum_{k=m+1}^\infty \frac{2}{\pi^2 k^2}\leq \frac{C}{ m}\, ,
\end{equation}
proving \eref{e:supnorm}.
\item
Let $Y_k$ be defined as in \eref{e:Yt}. Then, if ${\mathbb E}_0$ is the conditional expectation with respect to $X_0$
and $F_{X_k|X_0}$ is the conditional distribution function of $X_k$ given $X_0$,
\begin{align}
{\mathbb E}( \|{\mathbb E}_0(Y_k)\|_H)&={\mathbb E} \Big (
\Big(\int_0^1 (F_{X_k|X_0}(t)-F(t))^2 dt \Big)^{1/2}\Big) \nonumber \\
&\leq {\mathbb E}\Big( \sup_{t \in [0,1]}
|F_{X_k|X_0}(t)-F(t)|\Big):=\beta(k).
\label{beta}
\end{align}
Here  $\beta(k)$ is the weak $\beta$-mixing coefficient of the chain $(X_i)_{i \geq 0}$.
Starting from the computations of the paper by Dedecker,  Gou\"ezel and Merlev\`ede (2010),
we shall prove in the appendix that
\begin{equation}\label{maj}
\text{for any $\gamma \in ]0,1[$,} \quad
\beta (k) \leq \frac{C}{(k+1)^{(1-\gamma)/\gamma}} \, .
\end{equation}
\end{enumerate}

\subsection{Proof of Theorem \ref{mainth} for $\gamma=1/2$.} \label{gamma12}

Let $S_n$ be defined as in
\eref{e:Sn}. We shall prove that $Z_n=S_n/\sqrt{n\log (n)}$ satisfies
the points (i) and (ii) of Item 2 of Section \ref{Sec:gen}.

 We  first prove (i): for
any positive integer $m$,  $P_m(S_n/\sqrt{n \log (n)})$
converges in distribution in $H$  to
$P_m(V)$,  where $V=\{  g(t) Z, t \in [0,1]\}$ is the process
described in Item 1 of  Theorem \ref{mainth}.
For any ${\bf a}=(a_0, a_1, \ldots, a_m)$ in
${\mathbb R}^{m+1}$ and any
${\bf b}=(b_1, \ldots, b_m)$ in ${\mathbb R}^m$, let
\begin{equation}
  f_{{\bf a}, {\bf b}}(t)= a_0 + \sum_{k=1}^m a_k \cos(2 \pi k t) + \sum_{k=1}^m b_k
  \sin(2\pi k t) \label{f} \, .
\end{equation}
As noted in Section \ref{Sec:gen},
this is equivalent to prove that
\begin{equation}\label{e:f}
  \frac{1}{\sqrt{n \log (n)}} <f_{{\bf a}, {\bf b}}, S_n>
   \ \text{converges in distribution to}
  \ \sqrt{h(1/2)}\Big(\int_0^1 f_{{\bf a}, {\bf b}}(t)(1-F(t)) dt \Big) Z \, .
\end{equation}
Defining the function $u_{{\bf a}, {\bf b}}$
by
\begin{equation}
u_{{\bf a}, {\bf b}}(x)=\int_x^1 f_{{\bf a}, {\bf b}}(t) dt - \int_0^1
f_{{\bf a}, {\bf b}}(t) F(t) dt \, , \label{u}
\end{equation}
we obtain that
$$
<f_{{\bf a}, {\bf b}}, S_n>= \sum_{k=1}^n u_{{\bf a}, {\bf b}}(X_k) \, .
$$
Note that the function $u_{{\bf a}, {\bf b}}$ is Lipschitz
and that ${\mathbb E}(u_{{\bf a}, {\bf b}}(X_k))=0$. Hence, it follows from Gou\"ezel (2004)
that
$$
  \frac{1}{\sqrt{n \log (n)}} \sum_{k=1}^n u_{{\bf a}, {\bf b}}(X_k) \ \text{converges in distribution to}
  \ \sqrt{h(1/2)}u_{{\bf a}, {\bf b}}(0) Z \, .
$$
Since $u_{{\bf a}, {\bf b}}(0)=\int_0^1 f_{{\bf a}, {\bf b}}(t)(1-F(t)) dt$, \eref{e:f} holds and hence point (i) is proved.

We now prove (ii):
for any $\varepsilon>0$,
$$
\lim_{m \rightarrow \infty} \limsup_{n \rightarrow \infty}
{\mathbb P} \Big( \frac{\| S_n-P_m(S_n) \|_H^2}{\sqrt{n \log(n)}}> \varepsilon\Big)=0 \, .
$$
Hence, (ii) follows from
$$
\lim_{m \rightarrow \infty} \limsup_{n \rightarrow \infty}
\frac{{\mathbb E} ( \| S_n-P_m(S_n) \|_H^2)}{n \log (n)}=0  \, .
$$
By stationarity
\begin{align}
{\mathbb E} ( \| S_n-P_m(S_n) \|_H^2) &=
n  {\mathbb E}(\| Y_0-P_m(Y_0) \|_H^2)
+ 2 \sum_{k=1}^{n-1} (n-k){\mathbb E}(<Y_0-P_m(Y_0), Y_k-P_m(Y_k)>) \nonumber \\
&\leq n \Big( {\mathbb E}(< Y_0-P_m(Y_0), Y_0> )
+ 2 \sum_{k=1}^{n-1} |{\mathbb E}(<Y_0-P_m(Y_0), Y_k>)| \Big)\, ,
\label{e:key}
\end{align}
since $Y_0-P_m(Y_0)$ is orthogonal to $P_m(Y_k)$ for any $k=0, \ldots, n$.

Taking the conditional expectation  with respect to $X_0$, it follows that
$$
|{\mathbb E}(<Y_0-P_m(Y_0), Y_k>)|\leq |{\mathbb E}(<Y_0-P_m(Y_0), {\mathbb E}_0(Y_k)>)|
\leq \|\| Y_0-P_m(Y_0) \|_H \|_\infty{\mathbb E}( \|{\mathbb E}_0(Y_k)\|_H) \, .
$$
Therefore \eqref{e:key} yields
\begin{align}\label{ineqSn2}
{\mathbb E} ( \| S_n-P_m(S_n) \|_H^2)
 &\leq n \|\| Y_0-P_m(Y_0) \|_H \|_\infty  \Big(  {\mathbb E}(\|Y_0\|_H)
+ 2 \sum_{k=1}^{n-1} {\mathbb E}( \|{\mathbb E}_0(Y_k)\|_H) \Big) \nonumber \\
&\leq
2n \|\| Y_0-P_m(Y_0) \|_H \|_\infty  \sum_{k=0}^{n-1} {\mathbb E}( \|{\mathbb E}_0(Y_k)\|_H) \, .
\end{align}
By \eqref{e:supnorm}, $\|\|Y_0-P_m(Y_0)\|_H \|_\infty =O(m^{-1/2})$. To evaluate
the sum, we use the inequalities \eqref{beta} and
\eqref{maj}.  Since $\gamma=1/2$, it follows that
$\beta(k)=O(k^{-1})$, so that
\begin{equation}\label{sumbeta}
  \sum_{k=1}^{n-1} {\mathbb E}( \|{\mathbb E}_0(Y_k)\|_H) \leq \sum_{k=1}^{n-1} \beta(k) \leq C \log(n) \, .
\end{equation}
Consequently
$$
\lim_{m \rightarrow \infty}
\limsup_{n \rightarrow \infty} \frac{{\mathbb E} ( \| S_n-P_m(S_n) \|_H^2)}{n \log (n)}  \leq
\lim_{m \rightarrow \infty} \limsup_{n \rightarrow \infty}
\frac{C n \log(n)}{{\sqrt m} n \log(n)}=0 \, ,
$$
and (ii) follows.

\subsection{Proof of Theorem \ref{mainth} for   $\gamma \in ]1/2, 1[$.}\label{Sec:gamma>1:2}
\label{gamma>1/2}
Let $S_n$ be defined as in
\eref{e:Sn}. We shall prove that $Z_n=S_n/n^\gamma$ satisfies
the points (i) and (ii) of Item 2 of Section \ref{Sec:gen}.

We  first prove (i): for
any positive integer $m$,  $P_m(S_n/n^\gamma)$
converges in distribution in $H$  to
$P_m(V)$,  where $V=\{  g(t) Z, t \in [0,1]\}$ is the process
described in Item 2 of  Theorem \ref{mainth}.
For any ${\bf a}=(a_0, a_1, \ldots, a_m)$ in
${\mathbb R}^{m+1}$ and any
${\bf b}=(b_1, \ldots, b_m)$ in ${\mathbb R}^m$,
define the functions $f_{{\bf a}, {\bf b}}$ and
$u_{{\bf a}, {\bf b}}$ as in (\ref{f}) and (\ref{u}) respectively.  As in Section \ref{gamma12},
it suffices to prove that
$$
  \frac{1}{n^\gamma} \sum_{k=1}^n u_{{\bf a}, {\bf b}}(X_k) \ \text{converges in distribution to}
  \ C_\gamma (h(1/2))^\gamma \Big(\int_0^1
   f_{{\bf a}, {\bf b}}(t)(1-F(t)) dt \Big) Z \, .
$$
where $C_\gamma$ and $Z$ are described in Item 2 of  Theorem \ref{mainth}.
Note that the function $u_{{\bf a}, {\bf b}}$ is Lipschitz
and that ${\mathbb E}(u_{{\bf a}, {\bf b}}(X_k))=0$. Hence, it follows from Theorem 1.3
in  Gou\"ezel (2004)
that
$$
  \frac{1}{n¨^\gamma} \sum_{k=1}^n
  u_{{\bf a}, {\bf b}}(X_k) \ \text{converges in distribution to}
  \ C_\gamma (h(1/2))^\gamma u_{{\bf a}, {\bf b}}(0) Z \, .
$$
Since $u_{{\bf a}, {\bf b}}(0)=\int_0^1
f_{{\bf a}, {\bf b}}(t)(1-F(t)) dt$, the point (i) follows.

We now prove (ii):
for any $\varepsilon>0$,
$$
\lim_{m \rightarrow \infty} \limsup_{n \rightarrow \infty}
{\mathbb P} \Big( \frac{\| S_n-P_m(S_n) \|_H^2}{n^\gamma}> \varepsilon\Big)=0  \, .
$$
We shall  apply Proposition \ref{prop:maximal} of the appendix to the random variables
 $Y_i - P_m(Y_i)$ and the $\sigma$-algebras ${\mathcal F}_i= \sigma (X_1, \ldots, X_i)$.
 For any $j\geq i+k$, since ${\mathcal F}_i \subset {\mathcal F}_{j-k}$, one has
 $$
 {\mathbb E}(\|{\mathbb E}(Y_{j}-P_m(Y_j)|{\mathcal F}_i)\|_H) \leq {\mathbb E}(\|{\mathbb E}(Y_{j}-P_m(Y_j)|{\mathcal F}_{j-k})\|_H)\, .
 $$
Combined with the  Markov property, this implies that the coefficient $\theta(k)$ defined in Proposition \ref{prop:maximal} of the appendix
is such that: for $k \in \{0, \ldots, n-1\}$,
\begin{align*}
\theta(k)&=  \max \Big \{
 {\mathbb E}(\|{\mathbb E}(Y_{j}-P_m(Y_j)|{\mathcal F}_i)\|_H), (i, j)  \in \{1, \ldots , n\}^2
 \ \text{such that} \ j\geq i+k  \Big \} \\
 &=
 \max \Big \{
 {\mathbb E}(\|{\mathbb E}(Y_{j}-P_m(Y_j)|{\mathcal F}_{j-k})\|_H),  j  \in \{k+1, \ldots , n\}
  \Big \} \\
 &=
 \max \Big \{
 {\mathbb E}(\|{\mathbb E}(Y_{j}-P_m(Y_j)|X_{j-k})\|_H),  j  \in \{k+1, \ldots , n\}
   \Big \} \, .
\end{align*}
Let ${\mathbb E}_0$ be the conditional expectation with respect to $X_0$. By stationarity it follows that
\begin{equation}\label{e:tetha-k}
\theta (k) =  {\mathbb E}( \|{\mathbb E}_0(Y_k-P_m(Y_k))\|_H) \leq {\mathbb E}(\|{\mathbb E}_0(Y_k)\|_H)\, .
\end{equation}
the last inequality being satisfied because
$\|{\mathbb E}_0(Y_k)\|_H^2= \|{\mathbb E}_0(Y_k-P_m(Y_k))\|_H^2 +
\|{\mathbb E}_0(P_m(Y_k))\|_H^2$ by orthogonality.
By (\ref{maj})
 it follows that
$$
\theta (k) \leq \|{\mathbb E}_0(Y_k)\|_H \leq \beta (k) \, .
$$
Let $M(m)=
\|\| Y_0-P_m(Y_0) \|_H \|_\infty $ as in \eref{e:Mm}.
Applying Proposition \ref{prop:maximal} of the appendix,
for any positive integer $q$  and $x \geq  q M(m)$, one has
\begin{align}\label{DM}
{\mathbb P}\Big ( \max_{1 \leq k \leq n} \|S_k-P_m(S_k)\|_H \geq 4x
\Big ) &\leq
\frac{n \beta(q)}{x} + \frac{2nM(m)}{x^2} \sum_{k=0}^{q-1} \beta (q) \, .
\end{align}
By  (\ref{maj}), we know that $\beta (k)\leq
C (k+1)^{(\gamma-1)/\gamma}$.
Hence, it follows from \eqref{DM} that, for $x \geq qM(m)$ and $\gamma \in ]1/2, 1[$,
$$
{\mathbb P}\Big ( \max_{1 \leq k \leq n} \|S_k-P_m(S_k)\|_H \geq 4x
\Big ) \leq C \Big (
\frac {n}{x q^{(1-\gamma)/\gamma}} + \frac{q^{(2 \gamma -1)/\gamma}n
M(m)}{x^2}
\Big) \, .
$$
Taking $q=[x/M(m)]$ when $x\geq M(m)$,  we finally obtain that
\begin{align}\label{DMbis}
{\mathbb P}\Big ( \max_{1 \leq k \leq n} \|S_k-P_m(S_k)\|_H \geq 4x
\Big ) \leq \frac{C n M(m)^{(1-\gamma)/\gamma}}{x^{1/\gamma}}{\bf 1}_{x \geq M(m)}
+{\bf 1}_{x < M(m)}
\end{align}
(we bound this probability by 1 when $x < M(m)$).
We now apply \eqref{DMbis} with $x= n^\gamma \varepsilon / 4$. In view of the definition
\eref{e:Mm} of $M(m)$,
it follows that, for $n$ large enough,
$$
{\mathbb P} \Big( \frac{\| S_n-P_m(S_n) \|_H^2)}{n^\gamma}> \varepsilon\Big)\leq \frac{C}{\varepsilon^{1/\gamma}} \Big(\|\| Y_0-P_m(Y_0) \|_H \|_\infty
\Big)^{(1-\gamma)/\gamma} \,  .
$$
In view of  \eqref{e:bproj},
$$
\lim_{m \rightarrow \infty} \limsup_{n \rightarrow \infty}
{\mathbb P} \Big( \frac{\| S_n-P_m(S_n) \|_H^2)}{n^\gamma}> \varepsilon\Big)
\leq \frac{C}{\varepsilon^{1/\gamma}} \lim_{m \rightarrow \infty}
\Big(\frac{1}{m}\Big)^{(1-\gamma)/\gamma}= 0 \, ,
$$
since  $\gamma < 1$.
The point (ii)  follows.

\subsection{Proof of Theorem \ref{as}.}
Applying Inequality \eqref{equ2law}, we get that, for any $x>0$,
$$
\nu \Big( \max_{1 \leq k \leq n} k\|G_k\|_H \geq x \Big) \leq
{\mathbb P}\Big(2\max_{1 \leq k \leq n} \|S_k\|_H \geq x \Big) \,  .
$$
Hence, the inequalities  (\ref{dev1}) and (\ref{dev2}) hold provided the same
inequalities
hold for $\|S_k\|_H$ instead of $k\|G_k\|_H$.

Let $\gamma=1/2$. Replacing $Y_0-P_m(Y_0)$ by $Y_0$, we note that the inequality (\ref{DM}) is valid for $S_n$
with $M= \|\|Y_0\|_H\|_\infty$.
By (\ref{maj}), we know that $ \beta (q)\leq
C (q+1)^{-1}$.
Hence, it follows from \eqref{DM} that, for $x \geq qM$,
$$
{\mathbb P}\Big ( \max_{1 \leq k \leq n} \|S_k\|_H \geq 4x
\Big ) \leq C \Big (
\frac {n}{x q }+ \frac{\log(q) n
M}{x^2}
\Big) \, .
$$
Taking $q=[\sqrt{n}(\log(n))^\alpha]$ for some $\alpha \in ]0,1/2[$,  we finally obtain that
\begin{align}\label{DM2}
{\mathbb P}\Big ( \max_{1 \leq k \leq n} \|S_k\|_H \geq 4x
\Big ) \leq C \Big (
\frac {\sqrt n}{x (\log(n))^\alpha}+ \frac{\log(n) n
M}{x^2}
\Big){\bf 1}_{x \geq M[\sqrt{n}(\log(n))^\alpha]}
+{\bf 1}_{x < M [\sqrt{n}(\log(n))^\alpha]}\, .
\end{align}
Let $a_n$ be any  sequence of numbers  such that $a_n\geq a$
for some $a>0$, and
$\sum_{n>0} n^{-1} a_n^{-2}< \infty$.
Taking $4x=\varepsilon \sqrt{n\log(n)} a_n$ in \eqref{DM2}, we get that, for $n$ large enough
$$
{\mathbb P}\Big ( \max_{1 \leq k \leq n} \frac{\|S_k\|_H}
{\sqrt{n\log(n)} a_n} \geq \varepsilon \Big)
\leq C \Big (
\frac {1}{\varepsilon (\log(n))^{(1+2\alpha)/2}a_n}+
\frac{
M}{\varepsilon^2 a_n^2}
\Big ) \, .
$$
Since, by Cauchy-Schwarz,
$$
\sum_{n=2}^\infty \frac{1}{n} \frac {1}{ (\log(n))^{(1+2\alpha)/2}a_n}
\leq \Big ( \sum_{n=2}^\infty \frac{1}{n (\log(n))^{(1+2\alpha)}}
\Big)^{1/2} \Big (\sum_{n=2}^\infty  \frac{1}{n a_n^2}
\Big)^{1/2} < \infty \, ,
$$
we infer that
$$
\sum_{n=2}^\infty \frac{1}{n}
{\mathbb P}\Big ( \max_{1 \leq k \leq n} \frac{\|S_k\|_H}
{\sqrt{n\log(n)} a_n} \geq \varepsilon
\Big ) < \infty \, ,
$$
and (\ref{dev1}) follows. Assume moreover that $a_n$ is non decreasing  and such that
$a_n \leq c a_{n/2}$ for some $c\geq 1$.
Let $N\geq 2$ be a positive integer, and let
 $n \in \{2^N+1, \ldots,  2^{N+1}\}$. Clearly
\begin{equation}\label{encadrement}
\max_{1 \leq k \leq 2^N} \frac{k\|G_k\|_H}
{ \sqrt{2^N\log(2^N)}a_{2^N}} \leq
2c\max_{1 \leq k \leq n} \frac{k\|G_k\|_H}
{\sqrt{n\log(n)} a_n}
\leq 4 c^2 \max_{1 \leq k \leq 2^{N+1}} \frac{k\|G_k\|_H}
{ \sqrt{2^{(N+1)}\log(2^{N+1})}a_{2^{N+1}}}\, .
\end{equation}
Using the first inequality of \eqref{encadrement},  it follows from
(\ref{dev1}) that
\begin{align*}
\sum_{N=2}^\infty
{\nu}\Big ( \max_{1 \leq k \leq 2^N} \frac{k\|G_k\|_H}
{ \sqrt{2^{N}\log(2^N)}a_{2^N}} \geq \varepsilon
\Big ) &\leq \sum_{N=2}^\infty \sum_{n=2^N+1}^{2^{N+1}} \frac 2 n
{\nu}\Big ( \max_{1 \leq k \leq n} \frac{k\|G_k\|_H}
{\sqrt{n\log(n)} a_n} \geq \frac{\varepsilon}{2c}
\Big )\\
& \leq
\sum_{n=5}^{\infty} \frac 2 n
{\nu}\Big ( \max_{1 \leq k \leq n} \frac{k\|G_k\|_H}
{\sqrt{n\log(n)} a_n} \geq \frac{\varepsilon}{2c}
\Big ) < \infty \, .
\end{align*}
By the direct part of the Borel-Cantelli Lemma, and bearing in mind that $\nu$ is equivalent to the Lebesgue measure, we infer
that
\begin{equation}\label{e:2nd}
\lim_{N \rightarrow \infty} \max_{1 \leq k \leq 2^N} \frac{k\|G_k\|_H}
{ \sqrt{2^{N}\log(2^N)}a_{2^N}} =0 \quad
\text{almost everywhere}\,  .
\end{equation}
Using \eref{e:2nd} and the second inequality of \eqref{encadrement}, we conclude that
$$
  \frac{\sqrt n}
{\sqrt{\log(n)} a_n} \|G_n\|_H \quad \text{converges almost everywhere to 0}\,  ,
$$
proving \eref{e:asG1}.

Let $\gamma \in ]1/2, 1[$.
The inequality (\ref{DMbis}) is valid for $S_n$
with $M= \|\|Y_0\|_H\|_\infty$ and gives
\begin{equation}\label{DMter}
{\mathbb P}\Big ( \max_{1 \leq k \leq n} \|S_k\|_H \geq 4x
\Big )  \leq \frac{C n M^{(1-\gamma)/\gamma}}{x^{1/\gamma}}
{\bf 1}_{x \geq M}
+{\bf 1}_{x < M} \, .
\end{equation}
Let $a_n$ be any  sequence of positive numbers such that
$\sum_{n>0} n^{-1} a_n^{-1/\gamma} < \infty $.
Taking $x=4\varepsilon n^\gamma  a_n$ in \eqref{DMter}, we get that
$$
\sum_{n=1}^\infty \frac{1}{n}
{\mathbb P}\Big ( \max_{1 \leq k \leq n} \frac{\|S_k\|_H}
{n^\gamma a_n} \geq \varepsilon
\Big ) < \infty \, ,
$$
and (\ref{dev2}) follows.
Assume moreover that $a_n$ is non decreasing  and such that
$a_n \leq c a_{n/2}$ for some $c\geq 1$.
Using the same arguments as for the case $\gamma=1/2$, we infer
from (\ref{dev2}) that
$$
\sum_{N=1}^\infty
{\nu}\Big ( \max_{1 \leq k \leq 2^N} \frac{k\|G_k\|_H}
{2^{N\gamma}a_{2^N}} \geq \varepsilon
\Big ) < \infty \, .
$$
and we conclude that
$$
  \frac{n^{1-\gamma}}
{ a_n} \|G_n\|_H \quad \text{converges almost everywhere to 0}\,  ,
$$
proving \eref{e:asG2}.

\section{Extensions of the main results to some functions
of $T^k$.} \label{Sec:new}
\setcounter{equation}{0}
Let $g$ be  a function from $[0,1]$ to a compact interval $[a,b]$.
In this subsection, we modify the notations of Section \ref{intro}
 as follows:
we denote now by $G_n$ the empirical process of the sequence
$(g(T), g(T^2), \ldots, g(T^n))$, that is
\begin{equation}\label{e:Gbis}
G_n(t)=\frac 1 n \sum_{k=1}^n
\Big({\mathbf 1}_{g(T^k)\leq t}-F(t)\Big), \quad t \in {\mathbb R}
\end{equation}
where now $F(t)=\nu({\mathbf 1}_{g\leq t})$. Let also
$
  H={\mathbb L}^2([a,b], dt)
$.

\begin{thm}\label{prop:gene}
Let $g$ be a monotonic and H\"older continuous function from $[0,1]$ to ${\mathbb R}$. Then
the conclusions of Theorems \ref{mainth} and \ref{as} apply to the
process defined by \eref{e:Gbis}.
\end{thm}

\noindent{\it Proof of Theorem \ref{prop:gene}.}
Without loss of generality, assume that $g$ is a function
from $[0,1]$ to $[0,1]$, so that $H={\mathbb L}^2([0,1], dt)$.
The proof of this proposition is almost the same as that of
Theorems \ref{mainth} and \ref{as}. Let us  check the main points.

Let $S_n$ be the $H$-valued random variable defined by
$$
S_n(t)= \sum_{i=1}^n ({\mathbf 1}_{g(X_i) \leq t}-F(t))\, ,
$$
where $(X_k)_{k \geq 0}$ is a stationary  Markov chain with
invariant measure $\nu$ and transition kernel $K$ defined in
(\ref{Perron}).

Let us start with the convergence in distribution.

The finite dimensional convergence (Point (i) of Item 2 of Section \ref{Sec:gen}) can be proved as at
the beginning of Sections \ref{gamma12} and \ref{gamma>1/2}.
If $f_{{\bf a}, {\bf b}}$ and $u_{{\bf a}, {\bf b}}$ are defined by (\ref{f}) and
(\ref{u}) respectively, we obtain that
$$
<f_{{\bf a}, {\bf b}}, S_n>= \sum_{k=1}^n u_{{\bf a}, {\bf b}}(g(X_k)) \, .
$$
As already noticed, the function $u_{{\bf a}, {\bf b}}$ is Lipschitz, and consequently
the function $u_{{\bf a}, {\bf b}}\circ g$ is H\"older continuous. The finite
dimensional convergence follows as in Sections \ref{gamma12} and \ref{gamma>1/2}, since  Gou\"ezel's results (2004)
apply to any H\"older function.

The tightness
(Point (ii) of Item 2 of Section \ref{Sec:gen})
can be proved exactly as in Sections \ref{gamma12} and \ref{gamma>1/2}
provided that (\ref{e:supnorm}) holds for $Y_0={\mathbf 1}_{g(X_0) \leq t}-F(t)$, and provided that
the new coefficient
\begin{equation}\label{newbeta}
\beta(k)={\mathbb E}\Big( \sup_{t \in [0,1]}
|F_{g(X_k)|X_0}(t)-F(t)|\Big)
\end{equation}
satisfies (\ref{maj}). The first point can be proved as in (\ref{e:bproj}): for some positive constant $C$,
$$
\| Y_0-P_m(Y_0) \|^2_H= \sum_{k=m+1}^\infty (a_k^2(g(X_0))
+ b_k^2(g(X_0))) \leq \sum_{k=m+1}^\infty \frac{2}{\pi^2 k^2}\leq \frac{C}{ m}\, ,
$$
the Fourier coefficients $a_k$ and $b_k$ being defined in Section \ref{Sec:gen}.
To prove the second point note that, since $g$ is monotonic, the set
$\{g(X_k) \leq t\}$ is of the form $\{ X_k \leq u \}$, or $\{ X_k < u \}$, or  $\{X_k \geq u \}$, or $\{X_k > u\}$, for
some $u \in [0,1]$. Hence
$$
|F_{g(X_k)|X_0}(t)-F(t)|
\leq
 \sup_{u \in [0,1]}
|F_{X_k|X_0}(u)-{\mathbb P}(X_k \leq u)|\, ,
$$
and consequently
\begin{equation}\label{e:comp}\beta(k)=
{\mathbb E}\Big( \sup_{t \in [0,1]}
|F_{g(X_k)|X_0}(t)-F(t)|\Big)
\leq
{\mathbb E}\Big( \sup_{u \in [0,1]}
|F_{X_k|X_0}(u)-{\mathbb P}(X_k \leq u)|\Big)\, .
\end{equation}
By Proposition \ref{prop:majbeta} of the appendix
\begin{equation}\label{e:betabis}
{\mathbb E}\Big( \sup_{u \in [0,1]}
|F_{X_k|X_0}(u)-{\mathbb P}(X_k \leq u)|\Big) \leq \frac{C}{(k+1)^{(1-\gamma)/\gamma}}\, ,
\end{equation}
for some positive constant $C$.
From \eref{e:comp} and \eref{e:betabis}, it follows that the coefficient $\beta(k)$ defined in \eref{newbeta}  satisfies (\ref{maj}).

For the almost sure behavior of $G_n$, the proof is exactly the same as that of
Theorem \ref{as}, since the coefficient $\beta(k)$ defined in (\ref{newbeta}) satisfies (\ref{maj}).

\section{Appendix}
In this section, $C$ is a positive constant which may vary from
line to line.
\setcounter{equation}{0}

\subsection{A maximal inequality in Hilbert spaces}

The following proposition is used in the proof of Theorem \ref{mainth}.
It is adapted from Proposition 4 in Dedecker and Merlev\`ede (2007).

\begin{prop}\label{prop:maximal}
Let $Y_1, Y_2, \ldots, Y_n$ be $n$ random variables with values in a separable Hilbert
space $H$ with norm $\|\, \cdot \|_H$, such that
 ${\mathbb P}(\|Y_k\|_H\leq M)=1$ and ${\mathbb E}(Y_k)=0$ for
any $k \in \{1, \ldots , n\}$. Let ${\mathcal F}_1, \ldots ,
{\mathcal F}_n$ be an increasing filtration such that $Y_k$ is ${\mathcal F}_k$
mesurable for any $k \in \{1, \ldots , n\}$.
Let $S_n=\sum_{k=1}^n Y_k$, and for
$k \in \{0, \ldots, n-1 \}$, let
$$
\theta(k)= \max \Big \{
 {\mathbb E}(\|{\mathbb E}(Y_{j}|{\mathcal F}_i)\|_H), (i, j)  \in \{1, \ldots , n\}^2
 \ \text{such that} \ j\geq i+k  \Big \} \, .
$$
Then, for any $q \in \{1, \ldots, n\}$, and any $x\geq qM$, the following inequality holds
$$
{\mathbb P} \Big ( \max_{1 \leq k \leq n} \|S_k\|_H \geq 4x\Big)
\leq  \frac{n \theta(q)}{x}{\bf 1}_{q<n} + \frac{2nM}{x^2} \sum_{k=0}^{q-1} \theta (k)  \, .
$$
\end{prop}
\noindent {\it Proof of Proposition \ref{prop:maximal}.}
Let $S_0=0$ and define the random variables $U_i$ by:
 $U_i=S_{iq}-S_{(i-1)q}$ for
$i \in \{ 1, \ldots, [n/q]\}$ and $U_{[n/q]+1}=S_n-S_{q[n/q]}$. By Proposition 4 in Dedecker and Merlev\`ede (2007), for any $x\geq Mq$,
\begin{align}\label{propFlo}
{\mathbb P} \Big ( \max_{1 \leq k \leq n} \|S_k\|_H \geq 4x\Big)
&\leq \frac{1}{x} \sum_{i=3}^{[n/q]+1}
{\mathbb E}( \|{\mathbb E}(U_i|{\mathcal F}_{(i-2)q})\|_H) +
\frac{1}{x^2} \sum_{i=1}^{[n/q]+1}{\mathbb E}(\|U_i-{\mathbb E}(U_i|{\mathcal F}_{(i-2)q})\|_H^2) \nonumber  \\
&\leq
 \frac{1}{x} \sum_{i=3}^{[n/q]+1}
{\mathbb E}( \|{\mathbb E}(U_i|{\mathcal F}_{(i-2)q})\|_H) +
\frac{1}{x^2} \sum_{i=1}^{[n/q]+1}{\mathbb E}(\|U_i\|_H^2) \, ,
\end{align}
the second inequality being satisfied because
$\|U_i\|_H^2=\|U_i-{\mathbb E}(U_i|{\mathcal F}_{(i-2)q})\|_H^2+
\|{\mathbb E}(U_i|{\mathcal F}_{(i-2)q})\|_H^2$ by orthogonality.

To handle the first term in (\ref{propFlo}), note that $\theta(k)$
decreases with $k$ and,
according to the definition of $\theta (k)$: for $i \in \{ 1, \ldots, [n/q]\}$,
\begin{equation}\label{B1}
{\mathbb E}(\|{\mathbb E}(U_i|{\mathcal F}_{(i-2)q})\|_H) \leq
\sum_{j=(i-1)q+1}^{iq}\|{\mathbb E}(Y_j|{\mathcal F}_{(i-2)q})\|_H \leq
\sum_{j=(i-1)q+1}^{iq} \theta(j-(i-2)q)
 \leq q \theta(q) \, ,
\end{equation}
and
\begin{equation}\label{B2}
{\mathbb E}(\|{\mathbb E}(U_{[n/q]+1}|{\mathcal F}_{({[n/q]-1})q})\|_H)
\leq
\sum_{j=q[n/q]+1}^{n}\|{\mathbb E}(Y_j|
{\mathcal F}_{({[n/q]-1})q})\|_H
\leq (n-q[n/q])\theta(q)\, .
\end{equation}
From \eqref{B1}
and \eqref{B2}, and taking into account that the sum from ${i=3}$ to $[n/q]+1$ is
$0$ if $q=n$,
 we infer that
\begin{equation}\label{B2bis}
\sum_{i=3}^{[n/q]+1}
{\mathbb E}( \|{\mathbb E}(U_i|{\mathcal F}_{(i-2)q})\|_H)
\leq n \theta(q) {\bf 1}_{q<n}  \, .
\end{equation}

To handle the second term in (\ref{propFlo}), we start from the
 basic equalities:
 for $i \in \{ 1, \ldots, [n/q]\}$,
\begin{equation}\label{e:basic1}
{\mathbb E}(\|U_i\|_H^2) = \sum_{j=(i-1)q +1}^{iq} {\mathbb E}(\|Y_j\|_H^2) +
2\sum_{j=(i-1)q +1}^{iq}  \sum_{\ell=(i-1)q +1}^{j-1} {\mathbb E}(<Y_j, Y_\ell>)\, ,
\end{equation}
and
\begin{equation}\label{e:basic2}
{\mathbb E}(\|U_{[n/q]+1}\|_H^2)=
\sum_{q[n/q] +1}^{n} {\mathbb E}(\|Y_j\|_H^2) +
2\sum_{j=q[n/q] +1}^{n}  \sum_{\ell=q[n/q] +1}^{j-1} {\mathbb E}(<Y_j, Y_\ell>)\, .
\end{equation}
Taking the conditional expectation of
$Y_j$ with respect to ${\mathcal F}_\ell$ and proceeding exactly as to prove  \eqref{ineqSn2}, we obtain from
\eref{e:basic1} and \eref{e:basic2}
that: for $i \in \{ 1, \ldots, [n/q]\}$,
\begin{align}\label{B3}
{\mathbb E}(\|U_i\|_H^2)
& \leq  2 \sum_{j=(i-1)q +1}^{iq}
\sum_{\ell=(i-1)q +1}^j \|Y_\ell \|_\infty
\| {\mathbb E}(Y_j|{\mathcal F}_{\ell})\|_H \nonumber \\
 &\leq  2M \sum_{j=(i-1)q +1}^{iq}
\sum_{\ell=(i-1)q +1}^j \theta (j-\ell)
\leq 2Mq \sum_{k=0}^{q-1}  \theta (k)\,   ,
\end{align}
and
\begin{equation}\label{B4}
\quad {\mathbb E}(\|U_{[n/q]+1}\|_H^2)
\leq 2M(n-q[n/q]) \sum_{k=0}^{n-q[n/q]-1}  \theta (k) \leq
2M (n-q[n/q])
 \sum_{k=0}^{q-1}
\theta (k)\, .
\end{equation}
From \eqref{B3}
and \eqref{B4},
we infer that
\begin{equation}\label{B5}
\sum_{i=1}^{[n/q]+1}{\mathbb E}(\|U_i\|_H^2)=
{\mathbb E}(\|U_{[n/q]+1}\|_H^2) +
\sum_{i=1}^{[n/q]}{\mathbb E}(\|U_i\|_H^2)
\leq 2nM\sum_{k=0}^{q-1} \theta (k)\, .
\end{equation}
Starting from \eqref{propFlo} and using  the upper bounds  \eqref{B2bis}
and \eqref{B5}, Proposition \ref{prop:maximal} follows.

\subsection{Proof of the upper bound \eqref{maj} on the coefficient $\beta(k)$.}
\begin{prop}\label{prop:majbeta}
Let $\gamma \in ]0,1[$, and let $(X_i)_{i \geq 0}$ be a stationary Markov
chain with transition kernel $K$ defined in \eref{Perron} and invariant measure $\nu$.
Then the coefficient $\beta(k)$ defined in \eref{beta} satisfies the upper bound
\eref{maj}.
\end{prop}
\begin{rem} In fact, the upper bound  \eref{maj} holds for
Markov chains associated with the class of
generalized Pomeau-Manneville  (GPM) maps  introduced  in Definition 1.1 of
Dedecker, Gou\"{e}zel and Merlev\`ede (2010).
\end{rem}

\noindent{\it Proof of Proposition \ref{prop:majbeta}.}
 The proof is a slight modification of the proof of Proposition 1.16
in Dedecker, Gou\"{e}zel and Merlev\`ede (2010) and is included here for the sake of completeness.

If $f$ is supported in $[0,1]$, let $V(f)$ be the variation of
the function $f$, given by
  \begin{equation}\label{def:var}
  \Var(f)=\sup_{x_0<\dots<x_N} \sum_{i=1}^N | f(x_{i+1})-f(x_i)|\, ,
  \end{equation}
where the $x_i$'s are real numbers (not necessarily in $[0,1]$).
Let also $\|df\|$ denotes the variation norm of
the signed measure $df$ on $[0,1]$, which is defined as in
\eref{def:var} with all the $x_i$'s  in $[0,1]$. Note that, since $f$
is supported on $[0,1]$,
\begin{equation}\label{e:Vdf}
  \|f\|_\infty \leq \Var(f) = \|df\| +|f(0)|+|f(1)| \, ,
\end{equation}
where $\|f\|_\infty=\sup \{|f(x)|, x \in [0,1]\}$.
Note also that, for any $x \in [0,1]$,
\begin{equation}\label{e:majdf}
 |f(x)- \nu(f)| = \Big|\int_0^1(F(t)-{\bf 1}_{x \leq t}) df(t)\Big|\leq \|df\|\, .
\end{equation}
In particular, it follows from \eref{e:Vdf} and \eref{e:majdf} that
\begin{equation}\label{e:Vdf2}
\Var(f-\nu(f)) =  \|df\| +  |f(0)- \nu(f)|+  |f(1)- \nu(f)|
\leq 3 \|df\| \, .
\end{equation}

Let $K$ be the transition kernel defined in \eref{Perron}.
Recall that an equivalent definition of the coefficient $\beta(k)$ defined
in \eref{beta} is
\begin{equation}\label{betabis}
   \beta(k) = \nu \Big( \sup_{f : \|df\| \leq 1} |K^k f - \nu(f) | \Big ) \,
\end{equation}
(cf. Lemma 1 in Dedecker and Prieur (2005)).
Recall also that one has the decomposition
 \begin{equation}
  \label{eq_somme}
  K^n f = \sum_{a+k+b=n} A_a (\I_{(z_1, 1]}) \cdot \nu( B_b f) +
  \sum_{a+k+b=n} A_a E_k B_b f + C_n f \, ,
  \end{equation}
where the operators $A_n$, $B_n$, $C_n$ and $E_n$ and
the sequence $(z_n)_{n \geq 0}$ are defined in
Section 3 of the paper by  Dedecker, Gou\"{e}zel and Merlev\`ede (2010).
In particular, it is proved in this paper that
\begin{equation}\label{EandB}
\Var(E_k f)\leq \frac{C} {k^{(1-\gamma)/\gamma}} \Var(f) \quad \text{and} \quad
\Var(B_n f) \leq \frac{C \Var(f)}{(n+1)^{1/\gamma}}\, .
\end{equation}

Following the proof of Proposition 1.16
in Dedecker, Gou\"{e}zel and Merlev\`ede (2010), one has that
\begin{equation}\label{Cn}
 |C_n(f)| \leq C\|f\|_\infty K^n \I_{[0,z_{n+1}]}\, ,
\end{equation}
and
\begin{equation}\label{1Cn}
\nu(K^n \I_{[0,z_{n+1}]})= \nu( [0,z_{n+1}]) \leq \frac{C}{ (n+1)^{(1- \gamma)/\gamma}} \, .
\end{equation}

We now turn to the term $\sum_{a+k+b=n}A_a E_k B_b f$ in
\eqref{eq_somme}. Following the proof of Proposition 1.16
in Dedecker, Gou\"{e}zel and Merlev\`ede (2010), for any bounded
function $g$,
\begin{equation}\label{1An}
 |A_n(g)| \leq C\|g\|_\infty K^n \I_{(z_1,1]\cap T^{-1}[0,z_{n}]}\, ,
\end{equation}
and
\begin{equation}\label{1Anbis}
\nu(K^n \I_{(z_1,1]\cap T^{-1}[0,z_{n}]})
= \nu((z_1,1]\cap T^{-1}[0,z_{n}]) \leq \frac{C}{ (n+1)^{1/\gamma}} \, .
\end{equation}
Using successively \eqref{1An}, \eqref{e:Vdf} and \eqref{EandB},  we obtain that
  \begin{align}
  \label{AEB}
  \Big| \sum_{a+k+b=n} A_a E_k B_b f  \Big|  &\leq
  C \sum_{a+k+b=n}\|E_k B_b f\|_{\infty} K^a \I_{(z_1,1]\cap T^{-1}[0,z_{a}]}
  \nonumber
  \\&
 \leq  C \sum_{a+k+b=n}\Var(B_b f)\frac { K^a \I_{(z_1,1]\cap T^{-1}[0,z_{a}]}}
  {(k+1)^{(1-\gamma)/\gamma}} \nonumber
  \\&
  \leq  C\Var(f) \sum_{a+k+b=n}\frac { K^a \I_{(z_1,1]\cap T^{-1}[0,z_{a}]}}
  {(k+1)^{(1-\gamma)/\gamma}(b+1)^{1/\gamma}} \, .
  \end{align}

We now turn to the term
$\sum_{a+k+b=n}A_a(\I_{(z_1,1]})\cdot \nu(B_b f)$ in
\eqref{eq_somme}. From the displayed inequality right before (3.13) in Dedecker, Gou\"{e}zel and Merlev\`ede (2010),
one has
\begin{equation}\label{controlcentre}
  \left|\sum_{b=0}^{n-a} \nu( B_b f ) \right| =
   \left| \sum_{b>n-a} \nu( B_b f ) \right|
  \leq \sum_{b>n-a} \Var(B_b f)
  \leq \sum_{b>n-a}  \frac{C\Var(f)}{(b+1)^{1/\gamma}}
  \leq  \frac{D\Var(f)}{(n+1-a)^{(1-\gamma)/\gamma}}\, .
  \end{equation}
From \eqref{controlcentre} and  \eqref{1An},
we obtain
  \begin{equation}\label{secondtermbis}
  \left|\sum_{a=0}^n A_a (\I_{(z_1,1]}) \cdot \left(\sum_{b=0}^{n-a}
  \nu( B_b f)\right)\right|
  \leq
  C \Var(f) \sum_{a=0}^n \frac{K^a \I_{(z_1,1]\cap T^{-1}[0,z_{a}]}}{(n+1-a)^{(1-\gamma)/\gamma}}\, .
  \end{equation}

From \eref{e:Vdf2}, $\Var(f-\nu (f))\leq 3 \|df\|$. Hence, it follows from
\eqref{eq_somme},
\eqref{Cn}, \eqref{AEB} and \eqref{secondtermbis} that
\begin{multline}\label{mainbound}
|K^n(f-\nu (f))| \\
\leq C \|df\| \Big (K^n \I_{[0,z_{n+1}]}+ \sum_{a=0}^n \frac{K^a \I_{(z_1,1]\cap T^{-1}[0,z_{a}]}}{(n+1-a)^{(1-\gamma)/\gamma}}
+ \sum_{a+k+b=n}\frac { K^a \I_{(z_1,1]\cap T^{-1}[0,z_{a}]}}
  {(k+1)^{(1-\gamma)/\gamma}(b+1)^{1/\gamma}} \Big)\, .
\end{multline}
From \eref{betabis}, \eqref{mainbound}, \eqref{1Cn} and \eqref{1Anbis}, it follows that
\begin{multline}
\beta(n) \leq
C \Big( \frac{1}{ (n+1)^{(1- \gamma)/\gamma}} +
\sum_{a=0}^n \frac{1}{(a+1)^{1/\gamma}(n+1-a)^{(1-\gamma)/\gamma}}
\\
 + \sum_{a+k+b=n}\frac {1}{(a+1)^{1/\gamma}
  (k+1)^{(1-\gamma)/\gamma}(b+1)^{1/\gamma}}
  \Big) \, .
\end{multline}
All the sums on right hand being of the same order (see Lemma 3.2 of Dedecker, Gou\"ezel and Merlev\`ede (2010), and its application at the beginning of  the proof of their Proposition 1.15), it follows that
$$
\beta(n) \leq  \frac{C}{ (n+1)^{(1- \gamma)/\gamma}}\, ,
$$
and the proof is complete.

\medskip
\noindent\textbf{Acknowledgements.} Herold Dehling was partially supported by the Collaborative Research Project  {\em  Statistical Modeling of Nonlinear Dynamic Processes} (SFB 823) of the German Research Foundation DFG. Murad S. Taqqu was partially supported by the NSF grants DMS-1007616 and
DMS-1309009 at Boston University.

\end{document}